\documentclass[12pt,psamsfonts]{article}
\usepackage{amssymb}
\newtheorem{theo}{Theorem}[section]
\newtheorem{remarkk}[theo]{Remark}
\newenvironment{rem}{\begin{remarkk}\rm}{\end{remarkk}}
\newtheorem{definition}[theo]{Definition}
\newenvironment{devery easyfi}{\begin{definition}\rm}{\end{definition}}
\newtheorem{prop}[theo] {Proposition}

\newtheorem{lemma}[theo]{Lemma}
\newtheorem{example}[theo]{Example}

\newcommand{\ZZ}{\ensuremath{\mathbb{Z}}}

\newcommand{\HH}{\ensuremath{\mathbb{H}}}
\newcommand{\PP}{\ensuremath{\mathbb{P}}}

\newcommand{\SSS}{\ensuremath{\mathcal{S}}}

\newcommand{\A}{\ensuremath{\mathcal{A}}}

\newcommand{\ra}{\ensuremath{\rightarrow}}

\def\Bbb{\bf}

\def\C{{\Bbb C}}

\begin{document}
\title{Some new surfaces  with $p_g = q = 0$. }
\author{I. C. Bauer, F. Catanese \\
\small Mathematisches Institut der Universit\"at Bayreuth \\
\small Universit\"atsstr. 30 \\
\small 95447 Bayreuth}
\date{\today}
\maketitle

\section{Introduction}

It is well known that an algebraic curve of genus zero is isomorphic  to the projective line. The
search for an analogous statement in the case of algebraic  surfaces led Max Noether to
conjecture that a smooth regular (i.e.,
$q(S) = 0$) algebraic surface with vanishing geometric genus ($p_g(S)  = 0$) should be a rational
surface. The first counterexample to this conjecture was provided by  F. Enriques (\cite{enrMS},
I), who introduced the so called Enriques surfaces by considering the  normalization of sextic
surfaces in 3-space double along the edges of a tetrahedron.
  Nowadays a large number of surfaces of general type with $p_g = q =  0$ is known, but the first
ones  were constructed in the thirties  by   L. Campedelli and L. Godeaux (cf. \cite{Cam},
\cite{god}: in their honour  minimal surfaces  of general type with  $K^2 = 1$  are called
numerical Godeaux surfaces, and those with
$K^2 = 2$ are called numerical  Campedelli surfaces).

In the seventies, after rediscoveries of these old examples, many new  ones were found through
the efforts of several authors (cf. \cite{bpv}, pages 234-237 and references therein). In 
particular, in the spirit of Godeaux' method to produce interesting surfaces as quotients $S = 
Z/G$ of simpler surfaces by the free action of a finite group $G$,  A. Beauville proposed a very 
simple construction  by taking as $Z$ the product $ X = C_1 \times C_2$ of two curves of 
respective genera  $g_1, g_2
\geq 2$, together with an action of a group $G$ of order $(g_1 - 1 ) (g_2 -1 )$ (this method
produces surfaces with
$K^2 = 8$). He also gave an explicit example as quotient of two Fermat curves
  (in section
$2$, we shall indeed show that his example leads to exactly two non  isomorphic surfaces).

In this paper we  will discuss Beauville's construction starting  from the bottom, i.e., as the
datum of two appropriate coverings of $\PP^1$ and address the problem  of classification of these
surfaces, which we are unable for time limits to achieve in this note.

The interest on this issue stems from the open problem that David  Mumford set forth at the
Montreal Conference in 1980 : "Can a computer classify all surfaces  of general type with $p_g
=0$? Our purpose is to show how complex this question is (and probably computers are needed even
if one asks a more restricted question).

First of all, all known surfaces of general type  with $p_g = q = 0,  K^2 = 8$ are quotients
$\HH \times \HH / \Gamma$ of the product of two upper half planes by a discrete cocompact group. 
Besides the cited examples, there are also quotients which are not related to products of curves, 
and were constructed long ago by Kuga and Shavel using quaternion algebras (cf. \cite{kug},
\cite{shav}).

It is still a difficult open question whether one can have a fake  quadric, i.e., a surface of
general type which is homeomorphic to $\PP^1 \times \PP^1$.

Studying the special case where $S$ is the quotient of a product of  two curves, we want to show
how huge is the number of components of the corresponding moduli  space, and how detailed and
subtle the classification is.

An important feature is also the question of rigidity: for some of  these surfaces the moduli
space consists of one or two points (cf. \cite{cat00}, \cite{cat03}), for others it has strictly
positive dimension, and in any case the construction yields  connected components of the moduli
space.

  Surfaces with
$p_g = q = 0$ were  also investigated  from other points of view. We  would like to mention
several articles by M. Mendes Lopes and R. Pardini ( \cite{PardDP}, 
\cite{MLP1}, \cite{MLP2}) where the authors study the problem  of describing and classifying the
failure of  birationality of the bicanonical map.
\\

  We will  here classify all smooth algebraic surfaces $S = C_1 \times  C_2 / G$, where
$C_1$, $C_2$ are as above curves of genus at least two and $G$ is a  finite abelian group acting
freely on $C_1 \times C_2$ by a product action, and yielding a quotient surface 
with $p_g = q = 0$. In this case
$K^2$ has to be equal to
$8$, and we will see that there are already several cases. Our first  main result is

\begin{theo} Let $S$ be a surface with $p_g = q = 0$ isogenous to a  higher product $C_1 \times
C_2 / G$ of unmixed type. If $G$ is abelian, then $G$ is one of the following groups:
$(\mathbb{Z} /2 \mathbb{Z})^3$, $(\mathbb{Z} /2 \mathbb{Z})^4$, 
$(\mathbb{Z} /3 \mathbb{Z})^2$,
$(\mathbb{Z} /5 \mathbb{Z})^2$. \\ Each of these groups really occur.
\end{theo}

We will then give a complete description of the connected components  of the moduli space that
arise from these surfaces. We remark again that for $G = (\mathbb{Z}  /5 \mathbb{Z})^2$, we get
two isolated points, i.e., these surfaces are rigid, but there are two  different ones. For the
other cases the group determines a positive dimensional irreducible connected component of the
moduli space. \\

In section $4$ we calculate explicitly (theorem \ref{torsion}) the torsion group $ T(S) = H_1(S,
\ZZ)$ of our surfaces with $G$ abelian: it turns out that in some cases  $T(S)$, which has a
natural surjection onto $G$, is strictly bigger than $G$, but it is exactly $G = (\ZZ / 5 \ZZ)^2 $
for the two Beauville surfaces. Whence, these are only distinguished by their fundamental group,
and not by the first homology group.

A classification of surfaces with $p_g = q = 0$ isogenous to a  product (i.e., releasing the
hypothesis that the group  be abelian) is possible, but it is quite  complicated and there are
many more cases as we will show in the last section, where we give a  list of examples of
surfaces with $p_g = q = 0$ isogenous to a higher product $C_1 \times  C_2 / G$ with non abelian
$G$, some of them already known, others new. We will postpone the  complete classification to a
forthcoming article (in the non abelian case a non trivial problem is  also the one of determining
the Hurwitz equivalence or inequivalence of certain systems of  generators of a finite group).

Quite similar is the case of surfaces isogenous to a product of curves and with
$ q = p_g = 1$. The classification of these surfaces is also related  to the determination of the
so called non standard case for the non birationality of the  bicanonical map. We refer the
reader for this topic to the forthcoming Ph.D. Thesis of  F. Polizzi.

\section{Basic invariants of surfaces isogenous to a product.}

  Let $S$ be a smooth connected algebraic surface over the complex  numbers.\\ First we will
recall the notion of surfaces isogenous to a higher product of  curves. By prop. 3.11 of
\cite{cat00} the following two properties 1) and 2) of a surface are  equivalent.

\begin{definition}  A surface $S$ is said to be {\em isogenous to a  higher  product} if and only
if, equivalently, either

1)  $S$ admits a finite unramified covering which is isomorphic to a  product of curves of genera
at  least two, or \\ 2) $S$ is a quotient  $S := (C_1 \times C_2) /  G$, where the $C_i$'s are
curves of genus at least two, and $G$ is a finite group acting freely on 
$Z:= (C_1 \times C_2)$.

We have two cases: the {\rm  mixed case} where the action of $G$  exchanges the two factors (and
then $C_1 , C_2$ are isomorphic), and the {\rm   unmixed case} where 
$G$ acts via a product action.
\end{definition}

We recall briefly the following results from \cite{cat00} (cf. also \cite{cat03}):

\begin{itemize}
\item
Let $S$ be isogenous to a product, and let $S'$ be another surface with the same fundamental
group as $S$ and such that $K^2_S = K^2_{S'}$ (equivalently, $\chi(S) = \chi (S')$ or $e(S)=
e(S')$): then $S'$ is orientedly diffeomorphic to $S$ and either $S'$ or its complex conjugate
surface $\bar{S'}$ belongs to an irreducible smooth family, yielding a connected component of the
moduli space of surfaces of general type.

\item
There is a unique minimal realization $S := (C_1 \times C_2) /  G$ (i.e., the genera $g_1, g_2$
of the two curves  $C_1 , C_2$ are minimal).  It follows that $G$, $g_1$, $g_2$ are invariants of
the fundamental group of $S$. 

\item 
The minimal realization provides an explicit realization of the above family as the datum
of two branched coverings $C_i \ra C'_i = C_i/ G$ whose topological type is completely determined
by the two orbifold exact group sequences
$ 1 \longrightarrow \pi_1 (C_i)  \longrightarrow \Pi (i) \longrightarrow G \longrightarrow 1,
$ obtained from the fundamental group exact sequence of the quotient map $ C_1 \times C_2 \ra S$

$$(**)  1 \longrightarrow \Pi_{g_1} \times \Pi_{g_2} \longrightarrow \pi_1 (S) \longrightarrow G
\longrightarrow 1$$
by moding out the normal subgroup  $\Pi_{g_{i+1}}$ .

\end{itemize}

We obtain an easier picture in the case where $q(S) = 0$, or, equivalently, $C'_1 \cong C'_2
\cong \PP^1$.

\begin{definition}
1) Let $G$ be a group. Then a {\bf spherical system of generators} of $G$ ( {\bf S.G.S. of G})
is an ordered sequence $\A = (a_1, \dots a_n)$ of generators of $G$ with the property that their
product $ a_1 \dots a_n = 1$.

2) If we choose $n$ points $P_1, \dots P_n \in \PP^1$, and a geometric basis  $\gamma_1, \dots
\gamma_n$ of $\pi_1( \PP^1 - \{P_1, \dots P_n \}) $ ( $\gamma_i$ is a simple counterclockwise loop
around
$P_i$, and they follow each other by counterclockwise ordering around the base point), then a
S.G.S. of
$G$ determines a surjective homo0morphism $\psi: \pi_1( \PP^1 - \{P_1, \dots P_n \})  \ra G$.

Now, the braid group of the sphere $\pi_0( Diff ( \PP^1 - \{P_1, \dots P_n \}))$ operates 
on such homomorphisms, and their orbits are called Hurwitz equivalence classes of spherical
systems of generators.
\end{definition}

With the above notation we obtain 

\begin{theo}
Let $S$ be a surface isogenous to a product, of unmixed type and with $q(S)=0$.
Then to $S$ we attach its finite group $G$ (up to isomorphism) and the equivalence classes of an
unordered pair of two S.G.S.'s $\A, \A'$ of $G$, under the equivalence relation generated by

1) Hurwitz equivalence for $\A$,

1') Hurwitz equivalence for $\A'$,

2) simultaneous conjugation for $\A, \A'$, i.e., for $\phi \in Aut(G)$, we let 
 $(\A = (a_1, \dots a_n),  \A' = (a'_1, \dots a'_n))$  be equivalent to
$(\phi (\A) = (\phi (a_1), \dots \phi ( a_n)) , \phi( \A') = (\phi( a'_1), \dots \phi  (a'_n)))$.

Then two surfaces $S$, $S"$ are deformation equivalent if and only if the corresponding
equivalence classes of pairs of S.G.S.'s of $G$ are the same. 
\end{theo}

{\em Proof.} If $S$, $S'$ are deformation equivalent, then they have an isomorphic fundamental
group exact sequence $(**)$. Whence, we get pairs of isomorphic orbifold exact sequences,
compatible with an identification of $G$ with a fixed group. 
Now, the orbifold exact sequences determine homomorphisms  $\psi_1: \pi_1( \C'_1 - \{P_1, \dots
P_n \})  \ra G$, $\psi_2: \pi_1( \C'_2 - \{P'_1, \dots
P'_m \})  \ra G$. One sees immediately that these pairs are defined up to equivalence (for
instance, 2) follows by the fact that a $G$-covering space is determined by the kernel of the
surjection of the fundamental group onto $G$, and not by the specific homomorphism).

Conversely, we see easily that if the equivalence classes are the same, then the surfaces are
deformation equivalent. 

\hfill Q.E.D.

\begin{rem}
Observe that, if the group $G$ is abelian, then Hurwitz equivalence of $\A = (a_1, \dots a_n)$
is simply permutation equivalence of the sequence $ (a_1, \dots a_n)$.
\end{rem}

We shall assume throughout that we have a surface $S$ isogenous to a  higher product and that we
are in the unmixed case, thus we have a finite group $G$ acting on  two curves $C_1, C_2$  with
genera
$g_1, g_2 \geq 2$, and  acting freely  by the product action on
$Z := C_1 \times C_2$.

Since
$$ K_S^2 = 8 \chi (\mathcal{O}_S) = \frac{8 (g_1 - 1)(g_2 -1)}{|G|}
$$

the assumption $p_g(S) = q (S)= 0$ implies that $K_S^2 = 8$.

\begin{rem} We have the following elementary but crucial formulae: \\  1) $8(g_1-1)(g_2-1) =
K_{C_1 \times C_2}^2 = |G| \cdot K_S^2 = 8 \cdot |G|$, whence
$$ |G| = (g_1 - 1) (g_2 - 1).
$$  2) Since $q(S) = 0$ we have $C_i / G \cong \mathbb{P}^1$ for $i =  1,2$, so by the Hurwitz
formula we get:
$$ |G| = \frac{2}{-2 + \sum_j (1 - \frac{1}{m_j})} (g_i - 1),
$$ where $i = 1,2$ and $m_j$ is the  branching index of a branch  point $P_j$ of $C_i
\longrightarrow
\mathbb{P}^1$. In particular, in view of $1)$ it must hold:
$$
\frac{2}{-2 + \sum_j (1 - \frac{1}{m_j})} \in \mathbb{N}.
$$
\end{rem}

It is easy to see that the number of branch points of the two  coverings $C_i \longrightarrow
\mathbb{P}^1$ cannot be too high.

\begin{lemma} Let $S = C_1 \times C_2 / G$ be as above. Then the  number of branch points of each
covering $C_i \longrightarrow \mathbb{P}^1$ is at most eight.
\end{lemma}

{\em Proof.} Assume  e.g. that $C_1 \longrightarrow C_1 / G = 
\mathbb{P}^1$ has at least $9$ branch points. Then
$$ |G| \leq \frac{2}{-2 + \sum_{j=1}^9 (1 - \frac{1}{2})} (g_1 - 1) = 
\frac{4}{5} (g_1 - 1),
$$ contradicting $g_2 \geq 2$. Therefore we can have at most $8$ branch points.
  \hfill Q.E.D.

\section{The case: G abelian} We will assume from now on that $G$ is  a finite abelian group. In
this section we will show that the only abelian groups which give  rise to a surface isogenous to
a product $S = C_1 \times C_2 / G$, of unmixed type and with $p_g = q  = 0$, are $(\mathbb{Z} /2
\mathbb{Z})^3$,
$(\mathbb{Z} /2 \mathbb{Z})^4$, $(\mathbb{Z} /3 \mathbb{Z})^2$ and 
$(\mathbb{Z} /5
\mathbb{Z})^2$.\\

Our first step is to limit the order of the group $G$.

\begin{prop} Let $G$ be a finite abelian group and let $C$ be a  smooth algebraic curve of genus
$g \geq 2$ admitting an action of $G$ such that $C / G = 
\mathbb{P}^1$. We denote by $r$ the number of branch points of the morphism $C \longrightarrow C
/ G$. If 
$r \geq 4$ then
$$ |G| \leq 4 (g - 1)
$$  except for the case $r = 4$ and where the multiplicities of the  branch points are $(2, 2, 3,
3)$ (then $G = \mathbb{Z} / 6$).
\end{prop}

{\em Proof.} Recall that, by the Riemann existence theorem, giving a Galois  covering $ C \ra C /
G =
\mathbb{P}^1$, with branch points $P_1, \dots  P_r$ and branching  indices $m_1, \dots  m_r$ is
equivalent, in the case where $G$ is abelian, to giving
\begin{itemize}
\item Elements $a_1, \dots  a_r$ of $G$ of respective orders $m_1, \dots  m_r$ (here $a_i$ is the
image in $G$ of  a geometric loop around
$P_i$) such that
\item
$a_1 +  \dots a_r = 0$
\item
  $a_1, \dots  a_r$  generate $G$.
\end{itemize} 
Note that the elements  $a_1, \dots  a_r$  are unique up to ordering.

If $r \geq 5$, then $\sum (1 - \frac{1}{m_j}) \geq \frac{5}{2}$, whence
$$ |G| \leq \frac{2}{-2 + \frac{5}{2}} (g - 1) = 4 (g - 1).
$$ Therefore it remains to analyse the case $r = 4$. We assume that  the multiplicities are $m_1
\leq m_2 \leq m_3 \leq m_4$.
$(2, 2, 2, 2)$ is obviously not possible, since it contradicts $g 
\geq 2$. $(2, 2, 2, n)$, for
$n \geq 3$, is not possible, since $a_1 + a_2 +a_3 = - a_4$ has order 
$2$ contradicting the fact that $a_4$ has order $n$. Suppose now that $(m_1, \ldots , m_4) = (2, 
2, 3, n)$. Then
$$
\frac{2}{-2 + \sum_j (1 - \frac{1}{m_j})} = \frac{6n}{2n - 3} \leq 4,
$$ for $n \geq 6$. We remark that $n = 4$ or $5$ is not possible,  since $- a_4 = a_1 + a_2 +
a_3$ has order $3$ or $6$. Therefore the only possible case is $(2,  2, 3, 3)$. Here we have $a_1
+ a_3 = -(a_2 + a_4)$ has order $6$, whence $G = \mathbb{Z} / 6$ and 
$g = 2$. For the remaining cases $(2, 2, \geq 4, \geq 4)$, $(2, \geq 3, \geq 3, \geq3)$ and $( 
\geq 3, \geq 3, \geq 3, \geq 3)$ it is immediate that
$$
\frac{2}{-2 + \sum_j (1 - \frac{1}{m_j})} \leq 4.
$$
  \hfill Q.E.D. \\

\bigskip

Our second step is to show that the group $G$ cannot be cyclic:

\begin{prop}\label{cyclic} Let $S$ be a surface isogenous to a higher  product $C_1 \times C_2 /
G$ such that $q = 0$. Then $G$ cannot be cyclic.
\end{prop}

{\em Proof of prop. \ref{cyclic}.} Both maps $C_i \longrightarrow C_i  / G \cong \PP^1$ determine
the following situation: $G \cong \ZZ/d$ is generated by elements
  $a_1, \dots a_r$  of respective orders $m_1, \dots m_r$, respectively by elements
  $b_1, \dots b_s$  of respective orders $n_1, \dots n_s$. We claim that $G$ cannot act freely on
$C_1 \times C_2$.  In fact, the stabilizers of some  point in the first curve
$C_1$ are exactly the subgroups generated by some element $a_i$.  Since $G$ is cyclic, the union
$\SSS $ of the stabilizers is the set of elements whose order divides  some $m_i$. If $\SSS'$ is
the union of the stabilizers for the action on the  second curve $C_2$, we want
$ \SSS \cap \SSS' = \{ 0 \}$. This amounts to requiring that $ 
\forall i =1, \dots r, j= 1,
\dots s$, the integers
$m_i$ and $ n_j$ are relatively prime. The condition that the $a_i$'s  generate is however
equivalent to $d$ being the least common multiple of the $m_i$'s.  Since  $d$ is also the least
common multiple of the $n_j$'s, we obtain a contradiction.

  \hfill Q.E.D. \\

\bigskip

We proceed discussing the case $r = 3$, and we assume again that the  multiplicities are
$(m_1, m_2, m_3)$ with $m_1 \leq m_2 \leq m_3$.

\begin{rem} We observe that $g.c.d. (m_1, m_2) = 1$ implies that $m_3  = m_1 \cdot m_2$ and $G$ is
cyclic of order $m_3$, a posibility which was already excluded.
\end{rem}

We are now ready to prove the following:

\begin{prop} \label{order16} Let $S$ be a surface isogenous to a  higher product $C_1 \times C_2
/ G$ such that $p_g = q = 0$. Then either \\ 1) $g_1, g_2 \leq 5$,  i.e. $|G| = (g_1 - 1) (g_2 -
1) \leq 16$, \\ or \\ 2) $G = (\mathbb{Z} / 5 \mathbb{Z})^2$.
\end{prop}

Before proving the above proposition we will prove the following weaker form.

\begin{prop} Let $S$ be a surface isogenous to a higher product $C_1 
\times C_2 / G$ such that
$p_g = q = 0$. Then either \\ 1) $g_1, g_2 \leq 5$, i.e. $|G| = (g_1  - 1) (g_2 - 1) \leq 16$, \\
or \\ 2) for one of the two curves the datum $(m_1, m_2, m_3; G) $ of  branching orders plus
occuring group yields a priori only one of the  following possibilities:
\\ a)
$(2, 6, 6;\mathbb{Z} / 2 \mathbb{Z} \oplus \mathbb{Z} / 6 
\mathbb{Z})$, \\ b) $(2, 8, 8;\mathbb{Z} / 2
\mathbb{Z} \oplus \mathbb{Z} / 8 \mathbb{Z})$, \\ c) $(2, 12,  12;\mathbb{Z} / 2 \mathbb{Z}
\oplus \mathbb{Z} / 12 \mathbb{Z})$, \\ d) $(2, 20, 20;\mathbb{Z} / 2 
\mathbb{Z} \oplus
\mathbb{Z} / 20 \mathbb{Z})$, \\ e) $(3, 6, 6;\mathbb{Z} / 3 
\mathbb{Z} \oplus \mathbb{Z} / 6
\mathbb{Z})$, \\ f) $(4, 4, 4;\mathbb{Z} / 4 \mathbb{Z} \oplus 
\mathbb{Z} / 4 \mathbb{Z})$, \\ g) $(5, 5, 5;(\mathbb{Z} / 5 \mathbb{Z})^2)$.
\end{prop}

{\em Proof.} We have already seen that if $C \longrightarrow C / G  = 
\mathbb{P}^1$ has $\geq 4$ branch points, then
$$ |G| \leq 4 (g - 1)
$$  except for the case $r = 4$ and the multiplicities of the branch  points are $(2, 2, 3, 3)$
(then $G = \mathbb{Z} / 6$). But by prop. \ref{cyclic} we know that  this case cannot occur.
Therefore we can assume that $C_1 \longrightarrow C_1 / G$ has $r =3$  branch points. We write
again the multiplicities $(m_1, m_2, m_3)$ with $m_1 \leq m_2 \leq  m_3$. They correspond again
to elements $a_i \in G$ of order $m_i$, generating $G$ such that $a_1  + a_2 + a_3 = 0$. \\
$(2, 2, n)$ is not possible since then $-2 + \sum_j (1 - \frac{1}{m_j}) < 0$. \\

1) ${\bf m_1 = 2:}$ \\ Then $m_2 \geq 4$, since for $m_2 = 3$, we  must have $m_3 = 6$, whence
$\alpha := -2 + \sum_j (1 - \frac{1}{m_j}) = 0$.\\ If $m_2 = 4$, then $m_3 = 4$ and $\alpha = 0$,
which is not possible. \\
$m_2$ odd implies that $G$ is cyclic, so we can exclude all these  cases by prop. \ref{cyclic}. \\
If $m_2 = 6$, then $m_3 = 6$ and $G = \mathbb{Z} / 2 \mathbb{Z} 
\oplus \mathbb{Z} / 6
\mathbb{Z}$.
\\ If $m_2 = 8$, then $m_3 = 8$ and $G = \mathbb{Z} / 2 \mathbb{Z} 
\oplus \mathbb{Z} / 8
\mathbb{Z}$.
\\ For $m_2 = 10, 14, 16, 18$ we see that $m_3$ has to be equal to 
$m_2$, but then
$\frac{2}{\alpha}
\notin \mathbb{N}$. \\ If $m_2 = 12, 20$, then we are in the cases $c)$ resp. $d)$ of the  claim.
\\ We assume now that
$m_3 \geq m_2 \geq 22$. Then we have
$$ |G| \leq 2 ( -2 + \frac{1}{2} + \frac{21}{22} + \frac{21}{22})  ^{-1} (g_1 - 1) = \frac{44}{9}
(g_1 - 1).
$$ Therefore if $\frac{2}{\alpha} \in \mathbb{N}$, then 
$\frac{2}{\alpha} \leq 4$. \\

2) ${\bf m_1 = 3:}$ \\
$m_2 = 3$ implies $\alpha = 0$, whereas $m_2 = 4, 5$ imply that $G$  is cyclic. Therefore we can
assume $m_2 \geq 6$. \\
$m_2 = 6$ implies $m_3 = 6$ and $G = \mathbb{Z} / 3 \mathbb{Z} \oplus 
\mathbb{Z} / 6
\mathbb{Z}$. \\
$m_2 > 6$ implies that either $G$ is cyclic or $m_2 \geq 9$. If $m_3 
\geq m_2 \geq 9$, then
$$ |G| \leq 2 ( -2 + \frac{2}{3} + \frac{8}{9} + \frac{8}{9}) ^{-1}  (g_1 - 1) = \frac{9}{2} (g_1
- 1).
$$ Therefore if $\frac{2}{\alpha} \in \mathbb{N}$, then 
$\frac{2}{\alpha} \leq 4$. \\

3) ${\bf m_1 = 4:}$ \\
$m_2 = 4$ implies that $m_3 = 4$ and $G = \mathbb{Z} / 4 \mathbb{Z} 
\oplus \mathbb{Z} / 4
\mathbb{Z}$. \\
$m_2 = 5$ implies again that $G$ is cyclic, which is not possible;  therefore we can assume that
$m_3 \geq m_2 \geq 6$. Then
$$ |G| \leq 2 ( -2 + \frac{3}{4} + \frac{5}{6} + \frac{5}{6}) ^{-1}  (g_1 - 1) = \frac{24}{5}
(g_1 - 1).
$$ Therefore, if  $\frac{2}{\alpha} \in \mathbb{N}$, then 
$\frac{2}{\alpha} \leq 4$. \\

4) ${\bf m_1 = 5:}$ \\
$m_2 = 5$ implies that $m_3 = 5$ and $G = \mathbb{Z} / 5 \mathbb{Z} 
\oplus \mathbb{Z} / 5
\mathbb{Z}$. \\ Therefore we have $m_3 \geq m_2 \geq 6$. Then
$$ |G| \leq 2 ( -2 + \frac{4}{5} + \frac{5}{6} + \frac{5}{6}) ^{-1}  (g_1 - 1) = \frac{60}{14}
(g_1 - 1).
$$ Whence, if  $\frac{2}{\alpha} \in \mathbb{N}$, then 
$\frac{2}{\alpha} \leq 4$. \\

5) ${\bf m_1 \geq 6:}$ \\ In this case we have
$$ |G| \leq 2 ( -2 + \frac{15}{6}) ^{-1} (g_1 - 1) = 4 (g_1 - 1).
$$

Therefore we have proven our claim.

\hfill Q.E.D.

\bigskip

In order to prove proposition \ref{order16} we have now to exclude  the cases $2a) - 2f)$ of the
previous result. This will be done in the following lemma.

\begin{lemma} Let $S$ be a surface isogenous to a higher product $C_1 
\times C_2 / G$ such that
$p_g = q = 0$. Then $G$ cannot be one of the following groups: \\ a) 
$\mathbb{Z} / 2 \mathbb{Z}
\oplus \mathbb{Z} / 6 \mathbb{Z}$, \\ b) $\mathbb{Z} / 2 \mathbb{Z} 
\oplus \mathbb{Z} / 8
\mathbb{Z}$, \\ c) $\mathbb{Z} / 2 \mathbb{Z} \oplus \mathbb{Z} / 12 
\mathbb{Z}$, \\ d)
$\mathbb{Z} / 2 \mathbb{Z} \oplus \mathbb{Z} / 20 \mathbb{Z}$, \\ e) 
$\mathbb{Z} / 3 \mathbb{Z}
\oplus \mathbb{Z} / 6 \mathbb{Z}$, \\ f) $\mathbb{Z} / 4 \mathbb{Z} 
\oplus \mathbb{Z} / 4
\mathbb{Z}$.
\end{lemma}

{\em Proof.} \\ a) In this case the multiplicities of the branch  points for $C_1$ have to be
$(2, 6, 6)$. Then the union of the stabilizers is equal to $\{ (1,0),  (0,x), (1,5), (1,3), (1,1)
\}$. Therefore there are only $2$ elements of order $3$ left, and  they cannot generate $G$. So
there is no possibility that $G = \mathbb{Z} / 2 \mathbb{Z} \oplus 
\mathbb{Z} / 6 \mathbb{Z}$ acts freely on $C_1 \times C_2$. \\ The cases $b) - e)$ are excluded 
exactly in the same way. \\ f) Let $a_1$, $a_2$, $a_3$ generate $G = \mathbb{Z} / 4 \mathbb{Z} 
\oplus \mathbb{Z} / 4
\mathbb{Z}$. Then we can assume w.l.o.g. that $a_1$, $a_2$ is a 
$\mathbb{Z} / 4 \mathbb{Z}$ basis. But then, if $\Sigma$ resp $\Sigma '$ denotes the set of 
stabilizers of $C_1$ resp.
$C_2$, we have
$$
\sharp (\Sigma \cap ((\mathbb{Z} / 2)^2 - \{ 0 \} )) \geq 2,
$$   and the same for $\Sigma '$. In particular $\Sigma \cap \Sigma ' 
\neq \emptyset$. \hfill Q.E.D.
\\

This proves theorem \ref{order16}. \\ We are now ready to formulate  the main result of this
section.

\begin{theo}\label{possiblegroups} Let $S$ be a surface with $p_g = q  = 0$ isogenous to a higher
product $C_1 \times C_2 / G$. If $G$ is abelian, then $G$ is one of  the following groups:
$(\mathbb{Z} /2 \mathbb{Z})^3$, $(\mathbb{Z} /2 \mathbb{Z})^4$, 
$(\mathbb{Z} /3 \mathbb{Z})^2$,
$(\mathbb{Z} /5 \mathbb{Z})^2$.
\end{theo}

{\em Proof.} We know by our previous considerations that $G =  (\mathbb{Z} /5 \mathbb{Z})^2$ or
$|G| \leq 16$.\\
  Moreover, $G$ cannot be cyclic, whence $|G| \in \{ 4,8, 12, 16 \}$. 
\\ Obviously, $G = (\mathbb{Z} /2 \mathbb{Z})^2$ is not possible and this excludes the  case $|G|
= 4$. \\ If $|G| = 8$, then either $G = \mathbb{Z} /2 \mathbb{Z} \oplus \mathbb{Z} /4 
\mathbb{Z}$ or $G =(\mathbb{Z} /2 \mathbb{Z})^3$. Assume that $G = \mathbb{Z} /2 
\mathbb{Z} \oplus \mathbb{Z} /4
\mathbb{Z}$. Since $G$ is not generated by elements of order $2$,  there must be at least one
generator of order $4$ for each of the curves $C_1$ and $C_2$. But  there is exactly one non
trivial element, namely $(0,2)$, which is the double of any element of  order $4$. Hence the
stabilizers of the two curves cannot intersect trivially and therefore $G$  cannot act freely on
$C_1
\times C_2$.\\ If $|G| = 12$, then $G$ can only be $\mathbb{Z} /2 
\mathbb{Z} \oplus \mathbb{Z} /6 \mathbb{Z}$ and this case was excluded  before. \\ If $|G|
= 16$,  then $G$ is one of the following groups: $(\ZZ /2 \mathbb{Z})^2 \oplus \mathbb{Z} /4 
\mathbb{Z}$, $(\mathbb{Z} /4)^2$, $(\mathbb{Z} /2)^4$. $(\mathbb{Z} /4)^2$ was already  excluded
and $(\mathbb{Z} /2)^2
\mathbb{Z} \oplus \mathbb{Z} /4 \mathbb{Z}$ is excluded in the same  way as $(\ZZ /2
\mathbb{Z}) \oplus \mathbb{Z} /4 \mathbb{Z}$, since there is also  only one element which can be
the double of an element of order $4$. \hfill Q.E.D.

\section{The moduli of surfaces with $p_g = q = 0$ isogenous to a  higher product (with abelian
group).}

In this section we will show that the groups in theorem 
\ref{possiblegroups} really occur. More precisely, we will describe exactly the corresponding
moduli spaces.

\subsection{$G = (\mathbb{Z} /2 \mathbb{Z})^3$}

Since every element of $G$ has order $2$, we clearly need $r \geq 5$  branch points for each
covering $C_i \longrightarrow \mathbb{P}^1$. It is now easy to see  that $r_1 = 5$ and $r_2 =
6$.We denote by $\mathcal{S}_i$ the union of the stabilizers of the  covering $C_i
\longrightarrow \mathbb{P}^1$. Then, since  $\mathcal{S}_i$ contains a basis, it has cardinality
at least $3$. Since however $\mathcal{S}_1$, $\mathcal{S}_2$ are disjoint, and their union has
cardinality at most $7$, we see  that $\mathcal{S}_1$ must contain exactly
$4$ elements (since the sum of the five elements is zero). We may then  assume that
$$ (a_1, a_2, a_3, a_4, a_5) = (e_1, e_2, e_3, e_1, e_2 + e_3),
$$ where $e_1$, $e_2$, $e_3$ is a suitable $\mathbb{Z} /2 \mathbb{Z}$  - basis of $(\mathbb{Z} /2
\mathbb{Z})^3$. Then there is only one possibility (up to  permutation) left for $(b_1, b_2, b_3
b_4, b_5, b_6)$ , namely
$$ (b_1, b_2, b_3, b_4, b_5, b_6) = (e_1 + e_2, e_1 + e_3, e_1 + e_2  + e_3, e_1 + e_2, e_1 +
e_3, e_1 + e_2 + e_3).
$$

Therefore we have shown the following

\begin{theo} The surfaces with $p_g = 0$ isogenous to a product with  group $G = (\mathbb{Z} /2
\mathbb{Z})^3$ form an irreducible connected component of dimension 
$5$ in their moduli space.\end{theo}

\begin{rem} This result was already shown by R. Pardini in
\cite{PardDP}, where she classifies surfaces with$p_g = 0$, $K^2 = 8$, which
are double planes. In fact  the above surfaces are the only ones in our list
having non birational bicanonical map.
\end{rem}

\subsection{$G = (\mathbb{Z} /2 \mathbb{Z})^4$}

Again, since there are only elements of order $2$ in $G$ we see that the number of branch points
for each covering $C_i 
\longrightarrow
\mathbb{P}^1$ has to be at least $5$. But since $|G| \leq 4 (g_i -  1)$ for both curves, we see
that $r_1 = r_2 = 5$. For the first curve $C_1$ we can assume
$$ (a_1, a_2, a_3, a_4, a_5) = (e_1, e_2, e_3, e_4, e := e_1 + \ldots + e_4),
$$ where $e_1$, $e_2$, $e_3$, $e_4$ is a $\mathbb{Z} /2 \mathbb{Z}$ -  basis of $(\mathbb{Z} /2
\mathbb{Z})^4$. Then the problem reduces to finding $v_1, \ldots ,  v_5 \in G$ such that \\ 1)
$\sum_{i=1}^5 v_i = 0$; \\ 2) $rank \ (v_1, \ldots , v_5) = 4$; \\ 3) 
$v_i$ is of weight $w = 2$ or $3$ (since
$e$ is the only vector in $G$ of weight $4$ and the $e_i$'s are the  only vectors of weight
$w(e_i)$ equal to $1$ in $G$).

\begin{rem}
$\sum_{i=1}^5 v_i = 0$ implies that $\sum_{i=1}^5 w(v_i) \equiv  0(2)$. Therefore the number
$n_3$ of vectors of weight $3$ in $\{ v_1, \ldots , v_5 \}$ has to be even.
\end{rem}

\begin{lemma} Only the case $n_3 = 2$ is possible.
\end{lemma}

{\em Proof.} Since there are $4$ elements of weight $3$ in 
$(\mathbb{Z} /2 \mathbb{Z})^4$ we have to exclude the cases $n_3 = 4$, $n_3 =0$. Assume $n_3 =
4$. Then  w.l.o.g. $v_i = e + e_i$ for $i = 1, \ldots , 4$. But then $v_5 = \sum _{i=1}^4 v_i =
e$,  which contradicts $w(v_5) \in
\{2,3 \}$. \\ Assume that $n_3 =0$. But since
$$
\sum _{v : w(v) = 2} v = e,
$$ five of these $6$ vectors of weight $2$ can never have sum zero. 
\hfill Q.E.D.

\bigskip

Therefore without loss of generality we can assume that $v_1 = e +  e_1$, $v_2 = e + e_2$. Then
$v_3$, $v_4$, $v_5$ have all weight two and their sum is equal to 
$e_1 + e_2$. We observe that we cannot have: \\
$$ | \bigcup _{i=3} ^5 supp(v_i) | = 3,
$$ because this would imply $\sum _{i=3} ^5 v_i = 0$. Therefore we  can assume that $v_3 + v_4 =
e$ and then $v_5 = e_3 + e_4$. Then we have two possibilities for 
$v_3$ and $v_4$, namely
$$ v_3 = e_1 + e_3, \ \ \ v_4 = e_2 + e_4;
$$ or
$$ v_3 = e_1 + e_4, \ \ \ v_4 = e_2 + e_3.
$$ But these two possibilities give rise to isomorphic surfaces,  since they are equivalent by
the permutation of $v_1$ and $v_2$.\\ Therefore we have shown the following:

\begin{theo} The surfaces with $p_g = 0$ isogenous to a product with  group $G = (\mathbb{Z} /2
\mathbb{Z})^4$ form an irreducible connected component of dimension 
$4$ in their moduli space.
\end{theo}

\subsection{$G = (\mathbb{Z} /3 \mathbb{Z})^2$}

Examples of this type have already been given by Dolgachev in \cite{dolg}.

In this case $G - \{ 0 \} = (\mathbb{Z} /3 \mathbb{Z})^2 - \{ 0 \}$  has $8$ elements, whence the
union of the stabilizers of each covering has to consist of exactly  four elements, i.e.
$|\mathcal{S}_1| = |\mathcal{S}_2| = 4$. Moreover we know that the  number of branch points of
each covering $C_i \longrightarrow \mathbb{P}^1$ is $4$. Thus we have  up to permutation:
$$ (a_1, a_2, a_3, a_4) = (a, b, -a, -b),
$$ and
$$ (b_1, b_2, b_3, b_4) = (a', b', -a', -b'),
$$ where $a$, $b$ (resp. $a'$, $b'$) is a basis of $(\mathbb{Z} /3 
\mathbb{Z})^2$. Therefore we have shown the following

\begin{theo} The surfaces with $p_g = 0$ isogenous to a product with  group $G = (\mathbb{Z} /3
\mathbb{Z})^2$ form an irreducible connected component of dimension 
$2$ in their moduli space.
\end{theo}

\subsection{$G = (\mathbb{Z} /5 \mathbb{Z})^2$}

These surfaces are a particular case of examples that were introduced  by A. Beauville (cf.
\cite{Beauville}). \\ We see that we have for both coverings $C_i 
\longrightarrow \mathbb{P}^1$ $3$ branch points and the multiplicity is always $5$.
 In particular,
these  surfaces are rigid. We shall see that here we have two components of the moduli space,
i.e.  there are two non isomorphic Beauville surfaces.\\ In order to give a Beauville surface it
is  equivalent to give the following data:
\begin{itemize}
\item $a_1$, $a_2$, $a_3 \in G$ of order $5$
such that they generate $G$ and their 
sum is zero; \\
\item $b_1$, $b_2$, $b_3 \in G$ of order $5$
such that they generate $G$ and their 
sum is zero. \\
\end{itemize} Moreover they have to fulfill the following condition:
$$ (<a_1> \cup <a_2> \cup <a_3>) \cap (<b_1> \cup <b_2> \cup <b_3>) = \{ 0 \}.
$$

We denote the set of sixtuples $(a_1, a_2, a_3, b_1, b_2, b_3)$  satifying the above conditions
by $\tilde{\mathfrak{M}}$. On $\tilde{\mathfrak{M}}$ the group 
$\mathcal{G} := Gl(2,\mathbb{Z} /5 \mathbb{Z}) \times \mathfrak{S}_3 \times \mathfrak{S}_3$ acts
in  the natural way. We remark that
$| \mathcal{G} | = 24 \cdot 20 \cdot 6 \cdot 6$. \\

Up to a permutation of the $b_i$'s we can write every element of 
$\tilde{\mathfrak{M}}$ as
$(e_1, e_2, -(e_1 + e_2), \lambda (e_1 + 2e_2), \mu (3e_1 + 4e_2), 
\rho (e_1 + 4e_2))$. This is possible since $(\mathbb{Z} /5 \mathbb{Z})^2 - \{ <e_1>, <e_2>, <e_1 
+ e_2> \}$ determines exactly three stabilizer groups. Since $\lambda (e_1 + 2e_2) + \mu  (3e_1 +
4e_2) + \rho (e_1 + 4e_2) = 0$ we see immediately that $\lambda = \mu = \rho$, which  implies
that there are at most
$4$ different Beauville surfaces.

\begin{theo} There are exactly two non isomorphic surfaces with $p_g  = 0$ isogenous to a product
with group $G = (\mathbb{Z} /5 \mathbb{Z})^2$.
\end{theo}

{\em Proof.} Since the cardinality of $\tilde{\mathfrak{M}}$ is $24 
\cdot 20 \cdot 12 \cdot 2$, there are at leasts two  orbits of $\mathcal{G}$. But
obviously  the two elements $(e_1, e_2, -(e_1 + e_2), \lambda (e_1 + 2e_2), \lambda
(3e_1 + 4e_2), \lambda  (e_1 + 4e_2))$ and $(e_1, e_2, -(e_1 + e_2), - \lambda (e_1 +
2e_2), - \lambda (3e_1 + 4e_2), - 
\lambda (e_1 + 4e_2)) \in
\tilde{\mathfrak{M}}$ are equivalent under $(\varphi, (1 \ 2), (1 \  2)) \in \mathcal{G}$, where
$\varphi \in Gl(2,\mathbb{Z} /5 \mathbb{Z})$ is given by $\varphi  (e_1) = e_2$, $\varphi (e_2) =
e_1$. This proves the claim.

\hfill Q.E.D.

\section{$H_1(S, \mathbb{Z})$ for surfaces isogenous to a product with $G$ abelian}

In \cite{bpv}, p. 237, there is a list of examples of minimal surfaces of general type with $p_g
= q = 0$. While for $1 \leq K^2 \leq 6$ for each example the first homology group is given, in
the case $K^2 = 8, 9$ there is a question mark.

This motivated us to calculate $H_1 (S,
\mathbb{Z})$ for surfaces $S = C_1 \times C_2 /G$ isogenous to a higher product of unmixed type
with $G$ abelian.\\ Let's recall  again some facts from \cite{cat00}. Let $g_i$ be the genus of
the curve
$C_i$ and  denote by $\Pi_g$  the fundamental group of a compact Riemann surface of genus
$g$. Then we have the following exact sequence
$$(**) \ 1 \longrightarrow \Pi_{g_1} \times \Pi_{g_2} \longrightarrow \pi_1 (S) \longrightarrow G
\longrightarrow 1.
$$ Since we have assumed that $S$ is of unmixed type, $\Pi_{g_1}$ and $\Pi_{g_2}$ are both normal
subgroups of $\pi_1 (S)$. We define $\Pi (i+1) := \pi_1 (S) / \Pi_{g_i}$, where $i+1$ is
considered as element in $\mathbb{Z} /2 \mathbb{Z}$. Then we get two exact sequences
$$ 1 \longrightarrow \pi_1 (C_i)  \longrightarrow \Pi (i) \longrightarrow G \longrightarrow 1,
$$ which are exactly the orbifold fundamental group exact sequences of the coverings $C_i
\longrightarrow C_i / G =: C'_i$, in particular $\Pi (i) = \pi_1 ^{orb} (C'_i - B | m)$ is the
orbifold fundamental group of $C_i \longrightarrow C_i / G$ (for the definition and properties of
the orbifold fundamental group we refer again to \cite{cat00}). \\ We henceforth have an exact
sequence 
$$ 1 \longrightarrow \Pi_{g_1} \times \Pi_{g_2} \longrightarrow \Pi (1) \times \Pi (2)
\longrightarrow G \times G \longrightarrow 1,
$$ where $\pi_1 (S)$ is the inverse image of $\varphi : \Pi (1) \times \Pi (2) \longrightarrow G
\times G$ of $G$ diagonally embedded in $G \times G$. In particular, we have the
following exact sequence
$$ (*) \ \ \ \ 1 \longrightarrow \pi_1 (S) \longrightarrow \Pi (1) \times \Pi (2) \longrightarrow
G \longrightarrow 1,
$$ where $\Pi (1) \times \Pi (2) \longrightarrow G$ is the composition of $\varphi$ with the map
$G \times G \longrightarrow G$, $(a, b) \mapsto a - b$. \\

 We observe the following.

\begin{rem} Let 
$$ 1 \longrightarrow A \longrightarrow B \longrightarrow G \longrightarrow 1
$$  be an exact sequence of groups and assume $G$ to be abelian. Then the following sequence is
exact:
$$ 1 \longrightarrow A^{ab} \longrightarrow \overline{B} := B / [A,A] \longrightarrow G
\longrightarrow 1.
$$
and moreover the abelianization of $B$ equals the abelianization of $\overline{B} $.
\end{rem}

We apply this remark repeatedly: first to the exact sequence $(**)$, obtaining

$$ 1 \ra H_1 \times H_2 \ra  \overline{\pi_1 (S)} \ra G \ra 1 $$
which embeds into the exact sequence 
$$ 1 \ra H_1 \times H_2 \ra  \overline{\Pi (1)} \times \overline{\Pi (2)}  \ra G \times G \ra 1
.$$

It follows that an element  $(h_1, h_2) \in  H_1 \times H_2 $ is a commutator in
$\overline{\pi_1 (S)}$ if and only if $h_1$ is a commutator in  $\overline{\Pi (1)}$ and $h_2$ is
a commutator in  $\overline{\Pi (2)}$.

 Therefore, if we define $G_i$ as the abelianization of  $\Pi (i)$ , or
equivalently of 
$\overline{\Pi (i)}$,  the sequence 
$$ 0 \longrightarrow \pi_1 (S)^{ab} \longrightarrow G_1 \times G_2 \longrightarrow G
\longrightarrow 0.
$$ is exact and we have proven the following

\begin{prop} Let $S = C_1 \times C_2 /G$ be a surface isogenous to a higher product of unmixed
type with $G$ abelian. Furthermore,  denote by $G_i$ the abelianization of the orbifold
fundamental group of $C_i \longrightarrow C_i / G$. Then 
$$ H_1 (S, \mathbb{Z}) = ker (G_1 \times G_2
\longrightarrow G
\times G
\longrightarrow G),
$$ where the last map is obviously given by $(a, b) \mapsto a - b$.
\end{prop}

 In the rest of the paragraph we will use our classification result of the previous section in
order to calculate the torsion groups of all surfaces isogenous to a higher product with $G$
abelian and $p_g = 0$. 

\begin{theo}\label{torsion} Let $S$ be a surface with $p_g = q  = 0$ isogenous to a higher product
$C_1
\times C_2 / G$ (not of mixed type) and assume $G$ to be abelian. Then we get the following
values of
$H_1 (S, \mathbb{Z})$: \\

1) $H_1 (S, \mathbb{Z}) = (\mathbb{Z} / 2 \mathbb{Z})^6$ for $G = (\mathbb{Z} /2 \mathbb{Z})^3$,
\\

2) $H_1 (S, \mathbb{Z}) = (\mathbb{Z} /2 \mathbb{Z})^4$ for $G = (\mathbb{Z} /2 \mathbb{Z})^4$, \\

3) $H_1 (S, \mathbb{Z}) = (\mathbb{Z} /3 \mathbb{Z})^4$ for $G = (\mathbb{Z} /3 \mathbb{Z})^2$, \\

4) $H_1 (S, \mathbb{Z}) = (\mathbb{Z} /5 \mathbb{Z})^2$ for $G = (\mathbb{Z} /5 \mathbb{Z})^2$.
\end {theo}

 {\em Proof.} 1) In this case $C_1 \longrightarrow C_1 / G = \mathbb{P}^1$ has $5$ branch points
$p_1, \ldots p_5$ of multiplicities $(2, 2, 2, 2, 2)$, whence the orbifold fundamental group
$\Pi(1) = \pi_1 ^{orb} (\mathbb{P}^1 - \{ p_1, \ldots , p_5 \}; (2, 2, 2, 2, 2))$ equals $<a_1,
\ldots a_5 | \ a_1^2 = \ldots  = a_5^2 = a_1 \cdot \ldots \cdot a_5 = 1 > = <a_1,
\ldots a_4 | a_1^2 =
\ldots a_4^2 = (a_1 \cdot \ldots \cdot a_4 )^2 = 1 >$. Therefore $G_1 = \Pi(1)^{ab} = (\mathbb{Z}
/ 2 \mathbb{Z})^4$. Since $C_2 \longrightarrow C_2 / G = \mathbb{P}^1$ has $6$ branch points,
again each of multiplicitiy $2$, we see that $G_2 = (\mathbb{Z} / 2 \mathbb{Z})^5$. Therefore 
$$ H_1 (S, \mathbb{Z} ) = ker((\mathbb{Z} / 2 \mathbb{Z})^4 \oplus (\mathbb{Z} / 2 \mathbb{Z})^5
\longrightarrow (\mathbb{Z} / 2 \mathbb{Z})^3) = (\mathbb{Z} / 2 \mathbb{Z})^6.
$$

2) Here $C_i \longrightarrow C_i / G = \mathbb{P}^1$ has $5$ branch points $p_1, \ldots p_5$
resp. $q_1, \ldots q_5$ of multiplicities $(2, 2, 2, 2, 2)$, whence the orbifold fundamental
group $\Pi(i) = \pi_1 ^{orb} (\mathbb{P}^1 - \{ p_1, \ldots , p_5 \}; (2, 2, 2, 2, 2))$ (resp.
$\pi_1 ^{orb} (\mathbb{P}^1 - \{ q_1, \ldots , q_5 \}; (2, 2, 2, 2, 2))$  equals $<a_1, \ldots
a_5 : a_1^2 = \ldots a_5^2 = a_1 \cdot \ldots \cdot a_5 = 1 > = <a_1, \ldots a_4 : a_1^2 = \ldots
a_4^2 = (a_1 \cdot \ldots \cdot a_4 )^2 = 1 >$. Therefore $G_i = \Pi(1)^{ab} = (\mathbb{Z} / 2
\mathbb{Z})^4$.  Therefore 
$$ H_1 (S, \mathbb{Z} ) = ker((\mathbb{Z} / 2 \mathbb{Z})^4 \oplus (\mathbb{Z} / 2 \mathbb{Z})^4
\longrightarrow (\mathbb{Z} / 2 \mathbb{Z})^4) = (\mathbb{Z} / 2 \mathbb{Z})^4.
$$ 

3) Here $C_i \longrightarrow C_i / G = \mathbb{P}^1$ has $4$ branch points, all of multiplicity
$3$, and as above we see that $G_i = (\mathbb{Z} / 3 \mathbb{Z})^3)$, whence 
$$ H_1 (S, \mathbb{Z} ) = ker((\mathbb{Z} / 3 \mathbb{Z})^3 \oplus (\mathbb{Z} / 3 \mathbb{Z})^3
\longrightarrow (\mathbb{Z} / 3 \mathbb{Z})^2) = (\mathbb{Z} / 3 \mathbb{Z})^4.
$$ 

4) In this case $C_i \longrightarrow C_i / G = \mathbb{P}^1$ has $3$ branch points, all of
multiplicity $5$, and as before we see that $G_i = (\mathbb{Z} / 5 \mathbb{Z})^2)$, whence 
$$ H_1 (S, \mathbb{Z} ) = ker((\mathbb{Z} / 5 \mathbb{Z})^2 \oplus (\mathbb{Z} / 5 \mathbb{Z})^2
\longrightarrow (\mathbb{Z} / 5 \mathbb{Z})^2) = (\mathbb{Z} / 5 \mathbb{Z})^2.
$$  
\hfill  Q.E.D.

\section{Some new examples with $G$ non abelian}

We postpone the classification of surfaces $S$ with $p_g = 0$  isogenous to a higher product with
$G$ non abelian to a forthcoming article. \\ In the rest of the paper  we will give however some
new examples of surfaces isogenous to a product with non abelian  group. We remark that several
examples were already given by Mendes Lopes and Pardini (cf. 
\cite{PardDP}, \cite {MLP1}).
\\ We observe that in the non abelian case we cannot find such a low  upper bound for the
cardinality of the group $G$ (as in (2.4)), in fact we will exhibit examples 
of surfaces
$S$ with
$p_g = 0$ isogenous to a higher product with $G = \mathfrak{A}_5$ and $G = 
\mathfrak{S}_4$. The reason is that, in the case $r=3$, the branching indices here do  not need
to satisfy the condition that $m_3$ be a divisor of the least common multiple of $m_1$, $m_2$.

In the non abelian case, however, more restrictions come from the  condition that the
two stabilizer sets $\mathcal{S}_1$,$\mathcal{S}_2$ have an empty  intersection. In
fact, here
$\mathcal{S}_1$ is the union of the conjugacy classes of the cyclic  subgroups generated by
$a_1, \dots  a_n$. Therefore, knowledge of the conjugacy classes of 
$G$ will help to find examples, while knowledge of the branching indices plus Sylow's  theorems
help to show that some cases do not occur.
\\

\subsection{$G = \mathfrak{A}_5$}

Observe that in this case the group contains exactly three non  trivial conjugacy classes,
completely determined by the order of the elements in the class ( 
$m=2$ gives the class of the
$15$ double transpositions which form five Klein subgroups $\mathcal  K_i \cong (\ZZ / 2 \ZZ)^2$,
$m=3$ gives the conjugacy class of the $20$ three cycles, $m = 5$  yields the conjugacy class of
the $24$ five cycles).

It follows that for one of the two curves only one branching index can occur.

In this case the formulae of section $1$ read:
$$ |G| = 60 = (g_1 - 1) (g_2 - 1),
$$

$$ |G| = 60 = \frac{2}{-2 + \sum_j (1 - \frac{1}{m_j})} (g_i - 1).
$$ Denoting $\frac{2}{-2 + \sum_j (1 - \frac{1}{m_j})}$ by $\alpha_1$  resp. $\alpha_2$ we remark
that branching of pure type give the following values for $\alpha_i$:

$$ (2, \ldots , 2) =: 2^r \ \ \ \Longrightarrow \ \ \ \alpha_i = 
\frac{4}{r-4} \leq 4;
$$
$$ (3, \ldots , 3) =: 3^h \ \ \ \Longrightarrow \ \ \ \alpha_i = 
\frac{3}{h-3} \leq 3;
$$
$$ (5, \ldots , 5) =: 5^n \ \ \ \Longrightarrow \ \ \ \alpha_i = 
\frac{5}{2n-5} \leq 5.
$$

  Therefore we need ``mixed branching'' for at least one of the two curves.

Observe moreover that the integrality of $\alpha_i$ implies $r \in \{  5, 6 ,8\}$, $h \in \{ 4 , 6
\}$, $n \in \{ 3, 5\}$.

{\bf Example 1.} \\ For $C_1$ we take pure branching of type $3^4$,  i.e. $g_1 = 4$, $\alpha_1 =
3$, and for $C_2$ we take branching of type $(2, 5, 5)$, i.e. $g_2 =  21$, $\alpha_1 = 20$. \\
Since here obviously the union of the stabilizer subgroups for each  curve have trivial
intersection (remark that conjugating elements of order $2$, $3$, 
$5$, you get again elements of order $2$, $3$, $5$), the problem is reduced to finding elements 
$a_1$, $a_2$, $a_3$, $a_4$ of order three such that their product is $1$, generating 
$\mathfrak{A}_5$ and elements $b_1$,
$b_2$, $b_3$ of orders $(2, 5, 5)$, such that their product is $1$  generating $\mathfrak{A}_5$.
\\ 1) We set $a_1 = (123)$, $a_2 = (345)$, $a_3 = (432)$, $a_4 =  (215)$. It is obvious that
these are elements of order $3$ of $\mathfrak{A}_5$ and that their  product is $1$. Therefore it
remains to verify that $\mathfrak{A}_5$ is generated by these  elements. But we observe that $a_1
\cdot a_2 = (12345)$, which is an element of order $5$, and $a_1 
\cdot a_3 \cdot a_1 = (14)(23)$, which has order $2$. Therefore the subgroup generated by 
$a_1$, $a_2$, $a_3$, $a_4$ has order at least $30$. Since $\mathfrak{A}_5$ is simple it cannot 
have a subgroup of order
$30$, whence $a_1$, $a_2$, $a_3$, $a_4$ generate $\mathfrak{A}_5$. \\  2) We set $b_1 = (24)
(35)$, $b_2 = (21345)$, $b_3 = (12345)$. Obviously $b_1 \cdot b_2 
\cdot b_3 = 1$. In order to show that $b_1$, $b_2$, $b_3$ generate $\mathfrak{A}_5$ it suffices 
to find an element of order
$3$ in $<b_1, b_2, b_3>$. E.g. $b_3 \cdot b_1 \cdot b_3 = (152)$. \\  Therefore we have
constructed a surface $S = C_1 \times C_2 / \mathfrak{A}_5$, where 
$g(C_1) = 4$, $g(C_2) = 21$.
\\

In \cite{PardDP} the author gives another surface isogenous to a  product with group
$\mathfrak{A}_5$. This surfaces is obviously different to ours since  in her case $g(C_1) = 5$,
$g(C_2) = 16$. We will return to these examples later.

\bigskip

{\bf Example 2.} \\ For $C_1$ we take pure branching of type $5^3$,  i.e. $g_1 = 6$, $\alpha_1 =
5$, and for $C_2$ we take branching of type $(2, 2, 2, 3)$, i.e. $g_2  = 13$, $\alpha_1 = 12$. \\
Again the union of the stabilizer subgroups for each curve have  trivial intersection, hence we
have to find elements $a_1$, $a_2$, $a_3$ of order five such that  their product is $1$,
generating $\mathfrak{A}_5$ and elements $b_1$, $b_2$, $b_3$, $b_4$  of orders $(2, 2, 2, 3)$,
such that their product is $1$ generating $\mathfrak{A}_5$. \\ We set 
$a_1 = (12534)$, $a_2 = (12453)$, $a_3 = (12345)$; $b_1 = (12) (34)$, $b_2 = (24) (35)$, $b_3  =
(14) (35)$ and $b_4 = (234)$. It is now easy to see that these choices satisfy the required 
conditions and we obtain a new surface $S = C_1 \times C_2 / \mathfrak{A}_5$ with $g(C_1) = 6$, 
$g(C_2) = 13$.

\bigskip

{\bf Example 3.} \\ For $C_1$ we take pure branching of type $2^5$,  i.e. $g_1 = 5$, $\alpha_1 =
4$, and for $C_2$ we take branching of type $(3, 3, 5)$, i.e. $g_2 =  16$, $\alpha_1 = 15$. \\
Again the union of the stabilizer subgroups for each curve have  trivial intersection, hence we
have to find elements $a_1$, $a_2$, $a_3$, $a_4$, $a_5$ of order two  such that their product is
$1$, generating $\mathfrak{A}_5$ and elements $b_1$, $b_2$, $b_3$ of  orders $(3, 3, 5)$, such
that their product is $1$ generating $\mathfrak{A}_5$. \\ We set $a_1  = (12) (34))$, $a_2 = (13)
(24)$, $a_3 = (14) (23)$, $a_4 = (14) (25)$, $a_3 = (14) (25)$; $b_1  = (123)$, $b_2 = (345)$,
$b_3 = (54321)$. It is now easy to see that these choices satisfy the  required conditions and we
obtain a new surface $S = C_1 \times C_2 / \mathfrak{A}_5$ with 
$g(C_1) = 5$, $g(C_2) = 16$. These surfaces were already constructed by R. Pardini in
\cite{PardDP}. \\

\subsection{$G = \mathfrak{D}_4 \times \mathbb{Z} / 2 \mathbb{Z}$}

In order to avoid misunderstanding we note that for us 
$\mathfrak{D}_4$ is the group generated by $x$, $y$ with the relations $x^4 = y^2 = e$ and $yxy =
x^{-1}$.

Observe that in $\mathfrak{D}_4$ the centre consists of $\{ e,  x^2\}$, and there are three more
conjugacy classes, namely,  $\{ x,  x^{-1}\}$, $\{ y, y x^2\}$ and $\{x y, y x \}$.

\bigskip

{\bf Example 4.} \\ We will now rewrite an example which was already  constructed by R. Pardini
(cf. \cite{PardDP}) in our algebraic setting. For the curve $C_1$ we  take pure branching of type
$(2, 2, 2, 2, 2, 2)$,  whereas for $C_2$ we take  branching of type $(2, 2, 2,
4)$. Whence,  $g_1 = 9$, $ g_2 = 3$.\\
  We set $a_1 = (y, 0)$, $a_2 = (yx, 1)$, $a_3 = (yx^2, 0)$, $a_4 =  (yx, 1)$, $a_5 = (x^2, 1)$,
$a_6 = (x^2, 1)$. Then obviosly $a_1, \ldots a_6$ generate 
$\mathfrak{D}_4 \times \mathbb{Z} / 2
\mathbb{Z}$ and their product is $(e, 0)$. Furthermore for $C_2$ we  set $b_1 = (e, 1)$, $b_2 =
(y, 1)$, $b_3 = (xy, 0)$, $b_4 = (x, 0)$. Again these elements  generate $G$ and have trivial
product. We obtain thus a surface $S = C_1 \times C_2 / 
\mathfrak{D}_4 \times \mathbb{Z} / 2
\mathbb{Z}$ with $g(C_1) = 9$, $g(C_2) =3$.

\subsection{$\mathfrak{S}_4$}

Here there is only the following algebraic possibility of a surface 
$S = C_1 \times C_2 /
\mathfrak{S}_4$ with $p_g = 0$.

\bigskip

{\bf Example 5.} \\
  For the curve $C_1$ we take branching of type $(2, 2, 2, 2,2,2)$, whereas for
$C_2$ we take branching of type $(3, 4, 4)$. Whence,   
$g_1 = 13$, $g_2 = 3$.\\
  We set $a_1 = a_2 = (12)$, $a_3 = a_4  = (23)$, $a_5 = a_6 = (34)$.  Obviously
their product is
$1$ and they generate.

For
$C_2$, we set
$b_1 = (123)$,
$b_2 = (1234)$,
$b_3 = (1243)$. One immediately verifies that $b_1 b_2 b_3 = 1$. Moreover, if $H' :=
<b_1,
b_2, b_3>$, then $H'$  contains the transposition $(34)$
and acts double transitively, since it contains $b_1 = (123)$. Whence $H' =
\mathfrak{S}_4$. 

We obtain thus a surface $S = C_1 \times C_2 / 
\mathfrak{S}_4$ with $g(C_1) = 13$, $g(C_2) =3$. Again this example was already constructed by R. 
Pardini (cf. \cite{PardDP}).


\begin{thebibliography}{HorI-III}


\bibitem[BPV]{bpv} Barth, W., Peters, C., Van de Ven, A., {\em Compact complex surfaces.}
 Ergebnisse der Mathematik und ihrer Grenzgebiete (3). Springer-Verlag, Berlin,(1984).
 
\bibitem[Bea]{Beauville} A. Beauville, {\em   Surfaces alg\'ebriques complexes} Asterisque {\em
54} Soc. Math. France (1978).

\bibitem[Cam]{Cam} L. Campedelli, {\em Sopra alcuni piani doppi notevolicon curve di diramazione
del decimo ordine.} Atti Acad. Naz. Lincei {\em 15}, (1932), 536 - 542.

\bibitem[Cat00]{cat00} F. Catanese, {\em Fibred surfaces, varieties  isogenous to a product and
related moduli spaces.} Am. J. of Math. {\em 122}, 1 - 44 (2000).

\bibitem[Cat03]{cat03} F. Catanese,{\em Moduli spaces of surfaces and real structures.}  Ann. of
Math. {\em 158 } (2003), 539-554.

\bibitem[Dolg98]{dolg}{\em Surfaces with  $q = p_g = 0$.}
in 'C.I.M.E. 1977: Algebraic surfaces', Liguori, Napoli (1981), 247-266.

\bibitem[EnrMS]{enrMS} F. Enriques, {\em Memorie scelte di geometria, vol. I, II, III.} 
   Zanichelli, Bologna, (1956), 541 pp., (1959), 527 pp.,(1966), 456 pp. .

\bibitem[God]{god} L Godeaux, {\em  Les involutions cycliques appartenant \'a une surface
alg\'ebrique} Actual. Sci. Ind.,{\em  no. 270}, Hermann, Paris, (1935).

\bibitem[Ku]{kug} Kuga, M., {\em FAFA Note.}
 (1975).

\bibitem[MLP1]{MLP1} Mendes Lopes, M., Pardini, R. {\em The  bicanonical map of surfaces with
$p\sb g=0$ and $K\sp 2\geq 7$.} Bull. London Math. Soc. {\em 33}  (2001), no. 3, 265--274.

\bibitem[MLP2]{MLP2} M. Mendes Lopes, R. Pardini, {\em The  bicanonical map of surfaces with $p_g
= 0$ and $K^2 \geq 7$.} Bull. London Math. Soc. {\em 35} (2003), no.  3, 337--343.

\bibitem[Pa] {PardDP} Pardini, R. {\em The classification of double  planes of general type with
$K\sp 2=8$ and $p\sb g=0$.} J. Algebra {\em 259} (2003), no. 1, 95--118.



\bibitem[Sha]{shav} Shavel, I. H., {\em A class of algebraic  surfaces of general type
constructed from quaternion algebras.} Pacific J. Math. {\em 76},  (1978), no. 1, 221--245.
\end{thebibliography}
\end{document}